\input amstex
\frenchspacing
\documentstyle{amsppt}
\magnification=\magstep1
\NoBlackBoxes
\def\FF{\bold F}
\def\Gal{\mathop{\text{\rm Gal}}}
\def\coker{\mathop{\text{\rm coker}}}
\def\ker{\mathop{\text{\rm ker}}}
\def\lcm{\mathop{\text{\rm lcm}}}
\def\ord{\mathop{\text{\rm ord}}}
\def\mod{\bmod}
\def\congr{ \equiv }
\def\iso{ \cong }
\def\Que{\bold Q}
\def\Zee{\bold Z}
\newcount\refCount
\def\newref#1 {\advance\refCount by 1
\expandafter\edef\csname#1\endcsname{\the\refCount}}
\newref BA   
\newref BL   
\newref HO   
\newref LA   
\newref LE   
\newref MO   
\newref MSa  
\newref MSb  
\newref RM   
\newref PO   
\newref PSZ  
\newref SC   
\newref ST   
\newref STV  
\newref WA   
\topmatter
\title
A two-variable Artin conjecture
\endtitle
\author Pieter Moree and Peter Stevenhagen\endauthor
\address \hskip-\parindent Mathematisch Instituut,
Universiteit Leiden, P.O. Box 9512, 2300 RA Leiden, The Netherlands
\endaddress
\email moree\@math.leidenuniv.nl, psh\@math.leidenuniv.nl\endemail

\abstract
Let $a, b\in\Que^*$ be rational numbers that are multiplicatively
independent.
We study the natural density $\delta(a,b)$ of the set of primes $p$
for which the subgroup of~$\FF_p^*$ generated by $(a\mod p)$ contains
$(b\mod p)$.
It is shown that, under assumption of the generalized Riemann
hypothesis, the density $\delta(a,b)$ exists and equals a positive
rational multiple of the universal constant
$S=\prod_{p\text{ prime}} (1-p/(p^3-1))$.
An explicit value of $\delta(a,b)$ is given under mild conditions
on $a$ and $b$.
This extends and corrects earlier work of Stephens [\ST].

We also discuss the relevance of the result in the context of
second order linear recurrent sequences and some numerical
aspects of the determination of $\delta(a,b)$.
\endabstract
\endtopmatter

\document

\head 1. Introduction
\endhead

\noindent
Artin's original conjecture on primitive roots gives, for
each non-zero integer $a$, a conjectural value $\delta(a)$
of the density of the set
$$\{p \text{ prime}:\quad \langle a\mod p \rangle= \FF_p^*\}\tag1.1$$
inside the set of all primes.
It equals $\delta(a)=c_a\cdot\prod_{p\text{ prime}} (1-{1\over p(p-1)})$, where
$c_a$ is some explicit rational number that is positive whenever $a$
is not equal to $-1$ or to a square.
Artin's conjecture was proved by Hooley [\HO] under the assumption
of the generalized Riemann hypothesis.
Unconditionally, there is not a single value of~$a$ for which the set
of primes in (1.1) has been proved to be infinite (cf.~[\RM]).

In this paper we study the density $\delta(a, b)$ of the similarly
defined set
$$\{p \text{ prime}:\quad b\mod p\in\langle a\mod p \rangle\subset\FF_p^*\}
\tag1.2$$
inside the set of all primes.
In view of our application in Section 6, we allow $a$ and~$b$
to be non-zero {\it rational\/} numbers, and exclude the finitely
many primes dividing the numerators and denominators of $a$ and $b$
from consideration in (1.2).

If $a$ and $b$ satisfy a multiplicative relation $a^xb^y=1$ for
exponents $x, y\in\Zee$ that are not both equal to zero,
then one can prove unconditionally that $\delta(a,b)$ is a rational
number, and that it is positive in all but a few trivial cases [\STV].
We will therefore restrict to the case where $a$ and $b$ are
{\it multiplicatively independent\/} in $\Que^*$, i.e., no
non-trivial relation of the type above holds.
In this case, the following `two variable Artin conjecture'
has been proved unconditionally.
\proclaim
{Theorem 1}
Let $a, b\in\Que^*$ be multiplicatively independent.
Then the set of primes defined by $(1.2)$ is infinite.
\endproclaim\noindent
Theorem 1 is actually a special case of a theorem of P\'olya [\PO], and
we include its short and elementary proof in Section 6.
It does not show that the set (1.2) contains a subset of primes
of positive density.

We will mainly be concerned with the density of the set (1.2).
Here the basic result is the following.
\proclaim
{Theorem 2}
Let $a, b\in\Que^*$ be multiplicatively independent, and assume
the validity of the generalized Riemann hypothesis.
Then $\delta(a, b)$ exists and equals
$$\delta(a, b) = c_{a,b}\cdot\prod_{p\text{ prime}} (1-{p\over p^3-1})$$
for some positive rational constant $c_{a,b}$.
\endproclaim\noindent
The constant $c_{a,b}$ in the theorem depends on the degrees
of the number fields
$$F_{i,j}=\Que(\zeta_{ij}, a^{1/ij}, b^{1/i})\tag1.3$$
for $i, j\in\Zee_{>0}$.
Here $\zeta_{ij}$ denotes a primitive $ij$-th root of unity.
As an explicit formula for $c_{a,b}$ is rather cumbersome to write down,
we only compute it explicitly in the `generic case'
where the factor group $\Que^*/\langle-1, a, b\rangle$ is torsionfree.
This means that $\pm a^x b^y$ is not an $n$-th power in $\Que^*$ when
$x$ and $y$ are not both divisible by $n$.
\proclaim
{Theorem 3}
Let $a, b\in\Que^*$ be multiplicatively independent, and suppose that
the group $\Que^*/\langle-1, a, b\rangle$ is torsionfree.
Define $r(n)$ for $n\in\Zee_{\ne0}$ by
$$r(n)=\prod_{p|n\text{ prime}}
{-p^{4-3\,\ord_p(n)} \over p^3-p-1}.$$
Then the constant $c_{a,b}$ in Theorem $2$ has the value
$$c_{a,b}= 1+r(\lcm(2, \Delta(a)))+e(b)r(\Delta(b))+e(ab)r(\Delta(ab)).$$
Here $\Delta(x)\in\Zee$ denotes the discriminant of the quadratic field
$\Que(\sqrt x)$ for $x\in\Que^*$, and we put
$$e(x)=\cases
{3\over10}&\text{if $\Delta(x)$ is odd;}\\
1&\text{if $\Delta(x)$ is even.} \endcases
$$
\endproclaim\noindent
The universal constant $S=\prod_{p\text{ prime}} (1-p/(p^3-1))$
in Theorem 2, which is the analogue of Artin's constant
$A=\prod_{p\text{ prime}} (1-1/p(p-1))$ arising for the original
Artin conjecture, already occurs in Stephens' paper [\ST].
Our Theorem~2 occurs for positive coprime integers $a$ and $b$
that are not perfect powers as [\ST, Theorem 3], but the explicit
value for $c_{a,b}$ given there is involved and incorrect.
The analytic part of Stephens' proof, which we summarize in the next
section, is correct and generalizes in a rather straightforward way to our
more general situation.
His explicit evaluation of $c_{a,b}$ however, which is only carried through
in one out of the eight subcases distinguished in [\ST], is incorrect,
yielding an expression that is symmetric in $a$ and~$b$.
Our proof of Theorem 3 separates the elementary calculus of
double sums from the algebraic facts concerning the field degrees
$[F_{i,j}:\Que]$.

In the final two sections we address the relevance of our results
in the setting of second order recurrent sequences and
deal with some numerical aspects of the density $\delta(a, b)$.

\head 2. Results of Hooley and Stephens
\endhead

\noindent
The proof of the special case of Theorem~2 occurring in [\ST]
proceeds along the lines of Hooley's proof [\HO] of the original
Artin conjecture.

Artin's basic observation is that, for $a\in\Que^*$ arbitrary and
$p$ a prime number with $\ord_p(a)=0$,
the index $[\FF_p^*:\langle a\rangle]$ is divisible by $j$
if and only if $p$ splits completely in the splitting field
$F_j=\Que(\zeta_j, a^{1/j})$ of the polynomial $X^j-a$ over $\Que$.
By the Chebotarev density theorem, the set of these primes
has natural density $1/[F_j:\Que]$.
The primes $p$ for which $a$ is a primitive root modulo $p$ are
those primes that do {\it not\/} split completely in {\it any\/} of the
fields $F_j$ with $j>1$.
In fact, it suffices to require that $p$ does not split completely
in any field $F_j$ with $j$ {\it prime\/}.
For Artin's conjecture,
a standard inclusion-exclusion argument readily yields the heuristic
value
$$\delta(a)=\sum_{j=1}^\infty {\mu(j)\over [F_j:\Que]}.\tag2.1$$
The right hand side of (2.1) converges whenever $a$ is different from $\pm1$,
since in that case $[F_j:\Que]$ differs from its `approximate value'
$j\cdot \varphi(j)$, with $\varphi$ denoting Euler's $\varphi$-function,
by a factor that is easily bounded in terms of $a$.
In fact, we obtain an upper bound for the upper density of
the set (1.1) in this way.
In order to turn this heuristic argument into a proof, Hooley
employs estimates for the remainder term in the prime number theorem
for the fields $F_j$
that are currently only known to hold under assumption of the
generalized Riemann hypothesis.

In the situation of Theorem 2, one can find a conjectural value
for $\delta(a, b)$ in a similar way.
For each integer $i\ge1$, one considers the set of primes $p$
with $\ord_p(a)=\ord_p(b)=0$
for which the index $[\FF_p^*:\langle a\rangle]$ is {\it equal\/} to $i$
and the index $[\FF_p^*:\langle b\rangle]$ is {\it divisible\/}
by $i$.
These are the primes that split completely in the field
$F_{i,1}=\Que(\zeta_i, a^{1/i}, b^{1/i})$, but not in any of the
fields $F_{i,j}=\Que(\zeta_{ij}, a^{1/ij}, b^{1/i})$ with $j>1$.
As before, inclusion-exclusion yields a conjectural value
for the density $\delta_i(a,b)$ of this set of primes, and summing
over $i$ we get
$$\delta(a,b)=\sum_{i=1}^\infty \delta_i(a,b)\sum_{i=1}^\infty\sum_{j=1}^\infty {\mu(j)\over [F_{i,j}:\Que]}.\tag2.2$$
Note that $\delta_1(a,b)$ is nothing but the primitive root
density $\delta(a)$ from (2.1).

As in the case of (2.1), the right hand side of (2.2) converges if the
degrees $[F_{i,j}:\Que]$ are not too far from their `approximate values'
$i^2j\cdot\varphi(ij)$.
As we will see in Lemma 3.2, this is exactly what the hypothesis that $a$ and
$b$ be multiplicativily independent implies.

The proof of (2.2) by Stephens [\ST], which assumes the Riemann hypothesis
for each of the fields $F_{i,j}$, closely follows Hooley's argument in [\HO].
The restrictive hypotheses on integrality and coprimality
of $a$ and $b$ are not in any way essential.
The only requirement for the argument to work is that,
up to a factor that can be uniformly bounded from below by some
positive constant, $[F_{i,j}:\Que]$ behaves as $i^2j\cdot\varphi(ij)$.
We will show this in the next section, so there is no need for us to
elaborate any further on the proof of (2.2); we will merely be dealing
with degrees of radical extensions in order to prove Theorems 2 and 3.

The universal constant $S=\prod_{p\text{ prime}} (1-p/(p^3-1))$ is
the value of the right hand side of (2.2) obtained by substituting
$[F_{i,j}:\Que]=i^2j\cdot\varphi(ij)$, just like Artin's constant
$A=\prod_{p\text{ prime}} (1-1/p(p-1))$ is obtained from the right hand
side of (2.1) by putting $[F_j:\Que]=j\cdot\varphi(j)$.
The `correction factors' $c_a$ in Artin's conjecture and
$c_{a,b}$ in Theorem~2 measure the deviation of the field degrees
$[F_i:\Que]$ and $[F_{i,j}:\Que]$ from these values.
As in the case of Artin's conjecture, the basic problem is that
radical extensions of $\Que$ involving square roots are not in
general linearly disjoint from cyclotomic extensions.
In the case of Theorem 3, the situation is sufficiently simple
to allow an  expression for $c_{a,b}$ by a formula that
can be fitted on a single line.
As the proof of Theorem 2 will show, all other cases can in principle
be dealt with in a similar way.

For each prime $q$, the corresponding factor $1-1/q(q-1)$ in the
definition of $A$ has a well known interpretation:
it is, generically, the fraction of the primes $p$ for which
$\langle a \mod p\rangle$ generates a subgroup of $\FF_p^*$
of index not divisible by $q$.
In a similar way, the factor $1-q/(q^3-1)$ at a prime $q$ in the product
for~$S$ represents, generically,
the fraction of the primes $p$ for which the number of factors
$q$ in the index $[\FF_p^*:\langle a \mod p\rangle]$ is {\it at most\/}
the number of factors $q$ in $[\FF_p^*:\langle b \mod p\rangle]$.

As in the case of Artin's conjecture, we cannot prove unconditionally
that the set of primes in (1.2) has positive lower density.
The problem is that Chebotarev's density theorem only allows
us to simultaneously impose conditions of the type above
at primes $q$ for {\it finitely many\/} primes $q$.
As the fraction of primes that is eliminated by imposing a condition at $q$
is generically positive, one can prove unconditionally that there exists
a set of primes of positive density for which $b\mod p$ is {\it not\/}
contained in the subgroup $\langle a \mod p\rangle$ of $\FF_p^*$.
This was already noted by Schinzel [\SC].
For the set (1.2) itself, the best unconditional result available
is Theorem 1.

\head 3. Radical extensions
\endhead

\noindent
For Artin's original conjecture, one has to compute
the degree of $F_j=\Que(\zeta_j, a^{1/j})$, the number field obtained
by adjoining all $j$-th roots of an element $a\in\Que^*\setminus\{\pm1\}$
to $\Que$.
The result may be found in [\WA, Prop; 4.1].
The key observation is that if we take $a\ne\pm1$
such that $\Que^*/\langle-1, a\rangle$ is torsionfree,
then $a$ is a square in $\Que(\zeta_n)$ if and only if the discriminant
$\Delta(a)$ of $\Que(\sqrt a)$ divides $n$;
moreover, $a$ is not a $k$-th power
with $k>2$ in any cyclotomic extension of $\Que$.
If $j_1$ divides $j$ and $a$ is as above, then Kummer theory yields
$$[\Que(\zeta_j, a^{1/j_1}):\Que]=\cases
{1\over2}j_1\cdot\varphi(j)&\text{if $2$ divides $j_1$ and $\Delta(a)$ divides
             $j$;}\\
j_1\cdot\varphi(j)&\text{otherwise.} \endcases\tag3.1
$$
For arbitrary $a\in\Que^*\setminus\{\pm1\}$ and $j\in\Zee_{>0}$,
let $t=\gcd\{\ord_p(a): p$ prime$\}$ be the order
of the torsion subgroup of $\Que^*/\langle-1, a\rangle$, and take
$j_1=j/\gcd(j,t)$.
There are two cases.
If $a$ is a $t$-th power in $\Que^*$, say $a=a_1^t$,
we have $F_j=\Que(\zeta_j, a^{1/j})=\Que(\zeta_j, a_1^{1/j_1})$
and (3.1) can be applied directly.
If $a$ is not a $t$-th power, then $t$ is even and
$-a=a_1^t$ is a $t$-th power in $\Que^*$.
We have an extension
$$
\Que(\zeta_j, \zeta_{2j}a_1^{1/j_1})=F_j\subset F'_j=\Que(\zeta_{2j}, a_1^{1/j_1})
$$
of $F_j$ of degree at most 2 for which the degree
$[F'_j:\Que]$ is given by (3.1), and
one is left with the determination of $[F'_j:F_j]\in \{1,2\}$
as in~[\WA].
A somewhat subtle case distinction is necessitated by the peculiarity that
the element $-4$, which is not a square in $\Que^*$,
turns out to be equal to $(1+\zeta_4)^4$ in $\Que(\zeta_4)$.
Note that if we have $[F'_j:F_j]=2$, then $j$ and $t$ are even
and $2j_1$ divides $j$.
Thus, the `degree loss' $j\varphi(j)/[F_j:\Que]$ with respect
to the generic value $j\varphi(j)$ is always an integer dividing $2t$.
The factor 2 reflects the fact the kernel of the natural map
$$
\Que^*/{\Que^*}^k \to \Que(\zeta_k)^*/{\Que(\zeta_k)^*}^k
$$
is an abelian group annihilated by 2.
More precisely, it vanishes if $k$ is odd;
if $k$ is even it is generated by the elements $x^{k/2}{\Que^*}^k$
satisfying $\Delta(x)|k$ and, for $k\congr4\mod 8$, the
element $-2^{k/2}{\Que^*}^k$.

In our two variable setting, where we deal with the fields
$F_{i,j}=\Que(\zeta_{ij}, a^{1/ij}, b^{1/i})$, the statement in terms of the
map above easily leads to the following generalization.
\proclaim
{3.2. Proposition}
Let $a, b\in\Que^*$ be multiplicatively independent, and let $t$
be the order of the torsion subgroup of $\Que^*/\langle-1, a, b\rangle$.
Then for all $i, j\in\Zee_{>0}$, the quantity
$$f_{i, j}{i^2 j \varphi(ij)\over [\Que(\zeta_{ij}, a^{1/ij}, b^{1/i}):\Que]}$$
is a positive integer dividing $4t$.
In the torsionfree case $t=1$, it is equal to
\item{--}
the number of elements in $\{1,\Delta(a), \Delta(b), \Delta(ab)\}$
dividing $ij$ if $i$ is even;
\item{--}
the number of elements in $\{1,\lcm(2,\Delta(a))\}$
dividing $ij$ if $i$ is odd.
\endproclaim
\noindent{\bf Proof.}
Pick $i, j\in\Zee_{>0}$, and write $k=ij$.
Let $W_{i, j}\subset \Que^*$ be the subgroup generated
by $a$ and $b^j=b^{k/i}$, and $\overline W_{i, j}$ the image of
$W_{i, j}$ in $\Que^*/{\Que^*}^k$.
As the order of $\overline W_{i, j}$ divides $ik=i^2j$, we write
$i^2j=\#\overline W_{i, j}\cdot t_{i, j}$ with $t_{i, j}\in\Zee_{>0}$.
By Kummer theory, the degree of
$$F_{i,j}=\Que(\zeta_{ij}, a^{1/ij}, b^{1/i})=\Que(\zeta_k,\root k\of {W_{i, j}})$$
over $\Que(\zeta_k)$ equals $\#\psi[\overline W_{i, j}]$, with
$\psi: \overline W_{i, j}\to \Que(\zeta_k)^*/{\Que(\zeta_k)^*}^k$ the natural map.
We deduce that the `degree loss' $f_{i, j}$ for $F_{i, j}$ can be written as
$f_{i, j}=t_{i, j}\cdot \#\ker\psi$.
It is a decomposition of $f_{i, j}$ into a factor $t_{i, j}$
coming from `torsion in $\Que^*$'and a factor $\#\ker\psi$
measuring the additional torsion caused by the adjunction of $\zeta_k$.

As $\overline W_{i, j}$ is a finite abelian group on 2 generators and
$\ker\psi\subset \overline W_{i, j}$ is annihilated by 2,
it is clear that $\#\ker\psi$ divides 4.
In order to show that $t_{i, j}$ divides $t$, we let
$T\subset \Que^*$ be the inverse image of the torsion subgroup of
$\Que^*/\langle-1, a, b\rangle$ under the natural map
$\Que^*\to\Que^*/\langle-1, a, b\rangle$.
Then $T$ contains $V=\langle -1, a, b\rangle$ as a subgroup of index $t$, and
$\overline W_{i, j}$ is the subgroup of $T/T^k\subset \Que^*/{\Que^*}^k$
generated by $a$ and $b^j=b^{k/i}$.
The integer $t_{i, j}$ is the order of the kernel of the composed map
$${\langle a, b^j\rangle\over \langle a^k, b^k\rangle}{W_{i, j}\over W_{i, j}\cap V^k}\longrightarrow
V/V^k\longrightarrow T/T^k,$$
so it divides $\#\ker\left[V/V^k\to T/T^k\right]$.
As $V/V^k$ and $T/T^k$ are finite abelian groups of the same order,
$\#\ker f=\#\coker f= [T: VT^k]$ divides $[T : V]=t$.
This shows that $f_{i, j}$ divides $4t$.

Assume now that $\Que^*/\langle-1, a, b\rangle$ is torsionfree.
Then we have $T=V$ and $t_{i, j}=1$ in the argument above,
and $f_{i, j}=\#\ker\psi$.
Clearly, $\ker\psi$ vanishes if $k$ is odd.
If $k$ is even, the 2-torsion subgroup of $\overline W_{i, j}$ is generated by
$a^{k/2}$ if $i$ is odd and by $a^{k/2}$ and $b^{k/2}$ if $i$ is even.
As $\psi$ vanishes on the residue class of
$x^{k/2}\in \langle a^{k/2}, b^{k/2}\rangle$ in $\overline W_{i, j}$ if and only
if $\Delta(x)$ divides $k$, we arrive at
the value for $t_{i, j}$ given in the lemma. \hfill$\square$
\proclaim
{3.3. Corollary}
The integer $f_{i, j}$ in $3.2$ only depends on the greatest common divisors
$\gcd(i, 2t)$ and $\gcd(ij, 8st)$, where $s$ is the product of the primes
$p$ for which $\ord_p(a)$ and $\ord_p(b)$ are not both equal to $0$.
\endproclaim
\noindent
{\bf Proof.}
In the proof above, one needs $\gcd(k, t)=\gcd(ij, t)$ to determine the kernel
$\ker\left[V/V^k\longrightarrow T/T^k\right]$ and
$\gcd(i, t)$ to determine the order $t_{i, j}$ of the intersection of
$W_{i, j}=\langle a, b^{k/i}\rangle$ with this kernel.
As $\Delta(x)$ divides $4s$ for all $x\in V$, we can determine the
kernel of $V/V^k\to \Que(\zeta_k)^*/{\Que(\zeta_k)^*}^k$ if we
know $\gcd(k, 8s)=\gcd(ij, 8s)$.
Knowledge of the parity of $i$ enables us
to intersect this kernel with $W_{i, j}/(W_{i, j}\cap V^k)$, thus
yielding the second factor $\#\ker\psi$ in $t_{i,j}$.
\hfill$\square$

\head 4. Evaluation of the basic double sum
\endhead

\noindent
{}From (2.2), 3.2 and 3.3 it is clear that, in order to evaluate $\delta(a,b)$,
we need to evaluate for non-zero integers $m, n$ the double sum
$$
S_{m, n}=\sum_{{i=1\atop m|i}}^\infty \sum_{{j=1\atop mn|ij}}^\infty
{\mu(j)\over i^2j\varphi(ij)}.\tag4.1$$
This is a rather straightforward computation in elementary number theory
leading to the following result.
\proclaim
{4.2. Theorem}
For $m,n\in\Zee_{>0}$,
the value $S_{m, n}$ of the series in $(4.1)$ is
the rational multiple
$$S_{m, n}{S\over m^3n^3} \prod_{p|n} {-p^4\over p^3-p-1}
                \prod_{{p|m\atop p\nmid n}} {p^3+p^2\over p^3-p-1}$$
of the universal constant $S=\prod_{p\text{ prime}} (1-{p\over p^3-1})$
occurring in Theorem~$2$.
\endproclaim
\noindent{\bf Proof.}
We can sum over all $i\ge 1$ in (4.1) after substituting $mi$ for $i$.
Putting $ij=nd$ and summing over all $d\ge 1$ then yields
$$S_{m, n}={1\over m^2n^2}
\sum_{d=1}^\infty {1\over d^2 \varphi(mnd)} \sum_{j|nd} j\mu(j).$$
Writing $\widetilde x=\prod_{p|x} p$ for the largest squarefree divisor of $x$,
we have for any integer $x$
$$
\sum_{j|x} j\mu(j)=\sum_{j|\widetilde x} j\mu(j)\mu(\widetilde x)\sum_{j|\widetilde x} j\mu(\widetilde x/j)\mu(\widetilde x)\varphi(\widetilde x).$$
This enables us to write
$$S_{m, n}={1\over m^2n^2}\sum_{d=1}^\infty
{\mu(\widetilde {nd})\varphi(\widetilde {nd})\over d^2 \varphi(mnd)}{\mu(\widetilde n)\varphi(\widetilde n)\over m^2n^2\varphi(mn)}
\sum_{d=1}^\infty f(d),$$
where $f$ is the {\it multiplicative\/} function defined by
$$
f(d)={1\over d^2}\cdot{\varphi(mn)\over\varphi(mnd)}\cdot
{\mu(\widetilde {nd})\over \mu(\widetilde n)}\cdot
{\varphi(\widetilde {nd})\over \varphi(\widetilde n)}.
$$
As $\sum_{d=1}^\infty f(d)$ is absolutely convergent,
we can use the values
$$f(p^k)=\cases
-p^{1-3k}&\text{for $p\nmid mn$;}\\
p^{-3k}&\text{for $p|n$;}\\
-(p-1)p^{-3k}&\text{for $p|m$, $p\nmid n$} \endcases
$$
of $f$ on the prime powers $p^k$ with $k\ge1$ to obtain an Euler
product expansion
$$\eqalignno{
S_{m, n}&{\mu(\widetilde n)\varphi(\widetilde n)\over m^2n^2\varphi(mn)}
\prod_{p|n} {p^3\over p^3-1}
\prod_{p|m,\, p\nmid n} {p^3-p\over p^3-1}
\prod_{p\nmid mn} {p^3-p-1\over p^3-1}\cr
&{S\over m^3n^3} {mn\over\phi(mn)}
\prod_{p|n}{(1-p)p^3\over p^3-p-1}
\prod_{p|m,\, p\nmid n} {p^3-p\over p^3-p-1}\cr
&{S\over m^3n^3} \prod_{p|n} {-p^4\over p^3-p-1}
\prod_{p|m,\, p\nmid n}{p^3+p^2\over p^3-p-1}.&\square\cr
}$$

\proclaim
{4.3. Corollary}
Define $r(n)$ 
as in Theorem $3$.
Then we have
$S_{1,n}=r(n)S$ and
$$S'_{2,n}=\sum_{{i=1\atop 2|i}}^\infty \sum_{{j=1\atop n|ij}}^\infty
{\mu(j)\over i^2j\varphi(ij)}=\cases
 {3\over10}r(n)S&\text{if $n$ is odd;}\\
r(n)S&\text{if $4|n$.} \endcases
$$
\endproclaim
\noindent{\bf Proof.}
The first equality is immediate by taking $m=1$ in 4.2.

If $n$ is odd, we have $S'_{2,n}=S_{2, n}={3\over10}r(n)S$.
If $4$ divides $n$, the condition $2|i$ in the definition
of $S'_{2,n}$ is superfluous as $\mu(j)$ vanishes for $4|j$.
It therefore equals $S'_{2,n}=S_{1,n}=r(n)S$.
\hfill$\square$

\head 5. Proof of Theorems 2 and 3.
\endhead

\noindent
We now have everything at our disposal to prove Theorems 2 and 3.
\medskip\noindent
{\bf Proof of Theorem 2.\/}
We substitute the value of $[F_{i, j}:\Que]$ from Lemma 3.2
into the expression for $\delta(a,b)$ provided by (2.2) to obtain
$$\delta(a,b)=\sum_{i=1}^\infty\sum_{j=1}^\infty
                        f_{i, j}\ {\mu(j)\over i^2j\varphi(ij)}.\tag{5.1}$$
By 3.2, there are only finitely many values of $f_{i, j}$ that can occur,
namely the divisors of $4t$.
By 3.3, the value of $f_{i, j}$ only depends on the greatest common divisors
of $i$ and $ij$ with certain integers depending on $a$ and $b$.
It follows that the set of pairs $(i,j)$ for which
$f_{i, j}$ equals a given divisor of $4t$
can be characterized in terms of a finite number of divisibility criteria
on $i$ and $ij$.
This enables us to write $\delta(a,b)$ as an integral
linear combination of our basic sums $S_{m,n}$ for suitable values of $m$ and $n$.
As each of these sums is a rational multiple of $S$ by 4.2,
we conclude that $\delta(a,b)$ is itself a rational multiple of $S$.

It is not at all clear from the preceding argument that the resulting value for
$\delta(a, b)$ will always be positive.
{}From the expression $\delta(a, b)=\sum_{i=1}^\infty \delta_i(a,b)$
as a sum of non-negative terms in (2.2),
we see that it suffices to show that there is a value
of $i$ for which $\delta_i(a,b)$ is non-zero.
If $a$ is not a square, we can take $i=1$ as $\delta_1(a,b)=\delta(a)$ is then
positive by Hooley's result.
For arbitrary $a$, there can be many values of $i$ with $\delta_i(a,b)=0$.
In fact, for many $i$ one can construct $a$ that satisfy
$[{\bold F}_p^*:\langle a\rangle]\ne i$ for almost all $p$.
A list of such values of $i$ can be found in [\LE, (8.9)--(8.13)].
The smallest value that is not in the list is $i=24$, and we will
show that not only $\delta_{24}(a)$, but also
$\delta_{24}(a,b)$ is always positive.

We are interested in the primes $p$ for which
$[\FF_p:\langle a\rangle]$ equals 24 and $[\FF_p:\langle b\rangle]$
is divisible by 24.
Up to finitely many exceptions, these are the primes that
split completely in the field
$E=\Que(\zeta_{24}, \root 24 \of a, \root 24 \of b)$, but not in
any of its extensions $E_n=\Que(\zeta_{24n}, \root 24n\of a, \root 24 \of b)$
for $n>1$.
By the results of Lenstra [\LE, Theorem 4.1], the set of these primes
has positive density (under GRH) unless there is an obstruction
`at a finite level', i.e., an integer $h$ such that every
automorphism $\sigma$ of the extension $E\subset E_h$ is trivial on
$E_{n(\sigma)}$ for some divisor $n(\sigma)>1$ of $h$.
Thus, it suffices to show that for each squarefree integer $h$, there exists
$\sigma\in \Gal(E_h/\Que)$ satisfying the following 2 conditions:
\smallskip
\item{1.}
$\sigma$ is the identity on $E$;
\item{2.}
if $p$ is a prime dividing $h$ and $q$ is the largest power of $p$
dividing $24h$, then $\sigma$ is not the identity on $\Que(\zeta_q)$.
\smallskip\noindent
In order to construct such an automorphism, we observe that the
maximal subfield $E^{\text{ab}}\subset E$ that is abelian over $\Que$
has the property that $\Gal(E^{\text{ab}}/\Que)$ is an elementary abelian
2-group.
Assume without loss of generality that $6$ divides $h$, and
let $q$ be a prime power as in condition 2.
Take $\sigma_q$ to be any non-trivial automorphism of $\Que(\zeta_q)$ that
is a {\it square\/} in $\Gal(\Que(\zeta_q)/\Que)$.
As $q$ is not a divisor of 24, the group
$\Gal(\Que(\zeta_q)/\Que)\iso(\Zee/q\Zee)^*$ is not of exponent 2,
and such an element $\sigma_q$ exists.
Define $\sigma_0$ as the automorphism of $\Que(\zeta_{24h})$
with restrictions $\sigma|_{\Que(\zeta_q)}=\sigma_q$.
Then $\sigma_0$ is a square in $\Gal(\Que(\zeta_{24h})/\Que)$, so
it is the identity on $E\cap\Que(\zeta_{24h})\subset E^{\text{ab}}$.
This implies that there is a unique extension of $\sigma_0$ to $E(\zeta_{24h})$
that is the identity on $E$.
Any extension $\sigma$ of this automorphism of $E(\zeta_{24h})$ to $E_h$
now meets our requirements. \hfill$\square$
\medskip\noindent
{\bf Proof of Theorem 3.\/}
In the case where $a, b\in\Que^*$ are multiplicatively independent and
$\Que^*/\langle-1, a, b\rangle$ is torsionfree, we can use the
values of $f_{i, j}$ from Lemma 3.2
and rewrite (5.1) explicitly as an integral linear combination
of the type encountered in the proof of Theorem 2.
The sums of the type $S'_{2,n}$ from 4.3,
which single out the contribution to $S_{1,n}$ of the terms
with even $i$, can be used to obtain the compact expression
$$
\delta(a,b)S_{1,1}+S'_{2,\Delta(a)}+S'_{2,\Delta(b)}+S'_{2,\Delta(ab)}
+\left(S_{1,\lcm(2,\Delta(a))}- S'_{2,\lcm(2,\Delta(a))}\right).
$$
The sum of the three terms involving $\Delta(a)$ yields
$S_{1,\lcm(2,\Delta(a))}$, since we have
$S'_{2,\Delta(a)}=S'_{2,\lcm(2,\Delta(a))}$;
this is immediate if $\Delta(a)$ is even, and if $\Delta(a)$ is odd
the equality $S'_{2,\Delta(a)}=S'_{2,2\Delta(a)}$ follows directly from
the definition of the sum $S'_{2,n}$ in~4.3.
The explicit value of $c_{a, b}=\delta(a,b)/S$ now follows easily from
the two statements in~4.3.  \hfill$\square$

\head 6. Reformulation in terms of recurrent sequences
\endhead

\noindent
If we write the rational numbers $a$ and $b$ in Theorems 2 and 3
as $a=a_1/a_2$ and $b=b_1/b_2$, then we find that the set
of primes that divide some term of the integer sequence
$\{b_2a_1^n-b_1a_2^n\}_{n=0}^\infty$ has positive density.
This formulation in terms of integers is useful in proving
the unconditional result in Theorem 1.
The proof given below, which is entirely elementary,
is an easy extension of the argument of P\'olya occurring
in [\PSZ, Chapter 8, problem 107].
A generalization to integer sequences of the form
$\{\sum_{i=1}^k c_ia_i^n\}_{n=0}^\infty$ for arbitrary $k\ge2$
was given by P\'olya in [\PO].
\medskip\noindent
{\bf Proof of Theorem 1.}
Write $a=a_1/a_2$ and $b=b_1/b_2$ as above, and take
$\gcd(a_1, a_2)=\gcd(b_1,b_2)=1$.
We have $a_1\ne \pm a_2$ by the hypothesis that
$a, b\in\Que^*$ be multiplicativily independent,
so $|x_n|$ tends to infinity with $n$.
We need to show that the set $S$ of primes that
divide $x_n =b_2a_1^n-b_1a_2^n$ for some $n\ge0$ is infinite.

Suppose that $S$ is finite, and set $\ell=\varphi(|x_0|\cdot\prod_{p\in S} p)$.
Clearly, we have $\ell>0$.
We derive a contradiction by showing
that the sequence $\{x_{\ell n}\}_{n=0}^\infty$ is bounded.
As $S$ is finite, it suffices to show that $\ord_p(x_{\ell n})$
remains bounded as a function of $n$ for each $p\in S$.
Suppose that $p\in S$ is a prime that does {\it not\/} divide $a_1a_2$.
Then we have $\ord_p(x_{\ell n})=\ord_p(x_0)$ for all $n$ since
$$x_{\ell n}-x_0=b_2(a_1^{\ell n}-1)-b_1(a_2^{\ell n}-1)$$
is by the definition of $\ell$ divisible by $p^{\ord_p(x_0)+1}$.
Suppose that $p\in S$ is a prime dividing $a_1a_2$.
Then $p$ divides exactly one of $a_1$ and $a_2$, say $a_1$,
and we have $\ord_p(x_{\ell n})=\ord_p(b_1)$ for all sufficiently large $n$.
\hfill$\square$
\medskip\noindent
Integer sequences of the form $\{b_2a_1^n-b_1a_2^n\}_{n=0}^\infty$
are linear recurrent sequences of order 2.
They can be defined by the recursion
$x_{k+2}=(a_1+a_2)x_{k+1}-a_1a_2x_k$ for all $k\ge0$
and the initial values $x_0=b_2-b_1$ and $x_1=b_2a_1-b_1a_2$.

Much effort has been spent on the determination of the set of primes
dividing linear recurrent integer sequences, see~[\BL] and the references
given there.
In the case of second order sequences, our Theorems 1--3 lead to the
following result.
\proclaim
{Theorem 4}
Let $r, s\in\Que$ be rational numbers, and
${\Cal R}=\{x_n\}_{n=0}^\infty$ an integer sequence satisfying
the second order recursion $x_{k+2}=rx_{k+1}-sx_k$ for all $k\ge0$.
Suppose that $X^2-rX+s$ splits in $\Que[X]$, and that
${\Cal R}$ does not satisfy a first order recursion.
Then the set of primes that divide some term of $\Cal R$ is infinite;
if we assume the generalized Riemann hypothesis, it has positive density.
\endproclaim\noindent
{\bf Proof.}
Let $a_1, a_2\in\Que$ be the two roots of the polynomial $X^2-rX+s$.

Suppose first that we have $a_1\ne a_2$.
Then we have $x_n=b_2a_1^n-b_1a_2^n$ for certain $b_1, b_2\in\Que$.
As $X$ does not satisfy a first order recurrence, we have
$a_1a_2\ne0$ and $b_1b_2\ne0$.
After replacing, if necessary, the sequence $\{x_n\}_{n=0}^\infty$ by
$\{\lambda^n x_n\}_{n=0}^\infty$ for suitable $\lambda\in\Que^*$, we
may assume that $a_1$ and $a_2$ are coprime integers.
This only changes the set of primes that divide some term of $\Cal R$
by finitely many primes.
In a similar way, after replacing $\{x_n\}_{n=0}^\infty$ by
$\{\lambda x_n\}_{n=0}^\infty$ for suitable $\lambda\in\Que^*$,
we may assume that $b_1$ and $b_2$ are coprime integers.
If $a=a_1/a_2$ and $b=b_1/b_2$ are multiplicativily independent in
$\Que^*$, we are in the situation of Theorems 1 and 2, and we are done.
If $a$ and $b$ are multiplicativily dependent, then
the results on torsion sequences from [\STV] imply
unconditionally that the set of primes that divide some term of $\Cal R$
has positive density.

Suppose next that we are in the inseparable case $a_1=a_2$.
Then we have $x_n=(b_1+b_2n)a_1^n$ for certain $b_1, b_2\in\Que$, and
$b_2\ne0$ by assumption.
Now all primes that do not divide $b_2$ divide some term of $\Cal R$,
so we obtain a set of prime divisors of density 1.
\hfill$\square$
\medskip\noindent
The hypothesis that $\Cal R$ does not satisfy a first order recursion
in Theorem 4 is only there to exclude trivialities.
In order to remove the assumption that $X^2-rX+s$ splits in $\Que[X]$,
one needs to prove the analogues of our Theorems 1 and 2 for the set
$(1.2)$ in the case where $a$ and $b$ are elements of norm 1 in
a quadratic number field $K$ and $\FF_p$ is replaced by the ring of integers
of $K$ modulo the principal ideal~$(p)$.
It turns out that the inert primes lead to various complications.
The torsion case can be found in [\STV] and, for the special
case where $a$ comes from the fundamental unit of $K$, in [\MSa].
For a treatment of the non-torsion example proposed by Lagarias
[\LA, p. 451], we refer to [\MSb].

\head 7. Numerical data
\endhead

\noindent
Just like Artin's constant $A=\prod_{p\text{ prime}} (1-1/p(p-1))$, the
universal constant $S=\prod_{p\text{ prime}} (1-p/(p^3-1))$
in Theorem 2 is defined by a slowly converging product.
One can obtain good numerical approximations to $S$, such as the approximation
$$S\approx 0.57595~99688~92945~43964~31633~75492~49669~25065~13967~17649$$
up to 50 decimal digits, by expressing $-\log S$ as a rapidly converging series
involving the values $\zeta(d)$ of the Riemann zeta-function at arguments
$d\ge2$.
This is done for Artin's constant in [\BA], and for $S$ the expression
$$-\log S=-\log \zeta(3) +
\sum_{d=2}^{\infty}\sum_{k|d}\log \zeta(d)
{a_k\over d}\mu{d\overwithdelims()k}$$
is derived in [\MO].
Here $a_k$ is defined by its initial values $a_1=0$, $a_2=2$, $a_3=3$
and the recursion formula $a_{k+3}=a_{k+1}+a_k$ for $k\ge 1$.

It is not computationally feasible to determine the rational numbers
$c_{a, b}$ in Theorem 2 from numerical data.
In the torsionfree case occurring in Theorem 3,
the value of $c_{a, b}$ lies between the extremal values
$$c_{\min} = c_{2, 5}={9343\over 9520}\approx .981\quad\hbox{and}\qquad
  c_{\max} = c_{5, 3}={28001\over 27370}\approx 1.023.$$
For the 41535 primes contained in the interval $[7, 500\,000]$, one finds
that we have $\overline 5\in\langle2\mod p\rangle$ for 23498 primes
and $\overline3\in\langle5\mod p\rangle$ for 24429 primes.
When divided by $S$, these fractions are approximately equal to
$.9823$ and $1.0212$, respectively.
This shows that the deviations from $S$ are `numerically visible'
in a qualitative sense, but it also makes clear that one cannot determine
the fraction $c_{a, b}$ from such data.

\enddocument